\documentclass[11pt,reqno]{amsart}
\calclayout
\usepackage{amssymb}
\newif\ifaddpics
\ifaddpics\usepackage{graphicx}\fi
\makeatletter
\def\swappedhead@plain#1#2#3{%
  \thmnumber{\@upn{\mdseries #2}}\thmname{\@ifnotempty{#2}{. }#1}%
  \thmnote{ \textmd{\upshape(#3)}}}
\makeatother
\theoremstyle{plain}

\renewcommand{\Im}{\mbox{Im}}
\swapnumbers

\newtheorem{prop}[subsection]{Proposition}

\numberwithin{equation}{section}

\newcommand{\thismonth}{\ifcase\month\or
  January\or February\or March\or April\or May\or June\or
  July\or August\or September\or October\or November\or December\fi
  \space\number\year}
\makeatletter
\newcommand{\low}{\@ifnextchar^{}{^{\vphantom x}}}
\newcommand{\high}{\@ifnextchar_{}{_{\vphantom I}}}
\makeatother
\DeclareSymbolFont{script}{U}{eus}{m}{n}
\DeclareSymbolFontAlphabet{\mathscr}{script}
\DeclareMathSymbol{\EuWedge}{0}{script}{"5E}
\DeclareMathAlphabet{\mathrmsl}{OT1}{cmr}{m}{sl}
\newcommand{\bbsymb}[2]{\newcommand{#1}{{\mathbb{#2}}}}

\newcommand{\oper}[3][n]{\newcommand{#2}{\mathop
  {\mathrm{#3}\null}\ifx n#1\nolimits\else\limits\fi}}
\newcommand{\rsoper}[3][n]{\newcommand{#2}{\mathop
  {\mathrmsl{#3}\null}\ifx n#1\nolimits\else\limits\fi}}
\bbsymb\C{C} \bbsymb\F{F} \bbsymb\IH{H}\bbsymb\N{N} \bbsymb\IP{P}
\bbsymb\Q{Q} \bbsymb\R{R} \bbsymb\U{U} \bbsymb\V{V} \bbsymb\W{W}
\bbsymb\Z{Z}
 \bbsymb\E{E}
\newcommand{\ip}[1]{\langle #1 \rangle}
\renewcommand{\phi}{\varphi}

\newcommand{\hf}{{\mbox{\small$\frac{1}{2}$}}}

\newcommand{\hg}{\widehat{g}}

\newcommand{\rd}{\mbox{d}}

\renewcommand{\geq}{\geqslant} \renewcommand{\leq}{\leqslant}
\theoremstyle{definition}

\newcommand{\del}{\partial}                 
\newcommand{\db}{\overline\partial}       
\newcommand{\zb}{\overline z}
\newcommand{\wb}{\overline w}

\begin{document}
\title[Conic degeneration of PE metrics]{Some remarks on conic degeneration and bending \\ 
of Poincar\'e-Einstein metrics}
\author{Rafe Mazzeo and Michael Singer} 
\address{Department of Mathematics, Stanford University \\ Stanford, CA 94305 USA}
\email{mazzeo@math.stanford.edu}
\address{School of Mathematics\\ University of Edinburgh\\ King's Buildings, Mayfield Road\\ Edinburgh
EH9 3JZ\\ Scotland.}
\email{m.singer@ed.ac.uk}

\begin{abstract} Let $(M,g)$ be a compact K\"ahler-Einstein manifold with $c_1 > 0$. Denote by $K\to M$ 
the canonical line-bundle, with total space $X$, and $X_0$ the singular space obtained by blowing down 
$X$ along its zero section. We employ a construction by Page and Pope and discuss an interesting 
multi-parameter family of Poincar\'e--Einstein metrics on $X$. One $1$-parameter subfamily 
$\{g_t\}_{t>0}$ has the property that as $t\searrow 0$, $g_t$ converges to a PE metric $g_0$ on $X_0$ 
with conic singularity, while $t^{-1}g_t$ converges to a complete Ricci-flat K\"ahler metric $\hat{g}_0$ 
on $X$. Another $1$-parameter subfamily has an edge singularity along the zero section of $X$, with 
cone angle depending on the parameter, but has constant conformal infinity. These illustrate some 
unexpected features of the Poincar\'e-Einstein moduli space.
\end{abstract}

\maketitle
\section{Introduction}
The goal in this paper is to construct and examine some properties of
a family of Poincar\'e-Einstein (PE) (asymptotically hyperbolic)
metrics on the total space $X$ of a certain class of line bundles
$L\to M$, where $(M,\hg)$ is a compact K\"ahler-Einstein manifold with
$c_1 > 0$. Recall that a PE metric is a complete asymptotically
hyperbolic metric defined on the interior of a compact manifold with
boundary $X$, in a sense to be reviewed below, and to each such metric
$g$ is associated its conformal infinity ${\mathfrak c}(g)$, which is
a conformal class on $\del X$. Most of the PE metrics in the special
family considered here also have edge singularities along the zero
section of this total space, but a codimension one subfamily consists
of metrics which are smooth there. There are also interesting limits
of this family, for example Ricci flat metrics on $X$ (again possibly
with edge singularities on the zero section), obtained by letting
certain of the parameters go  to zero or infinity.  The ansatz leading to
this family is due to Page and Pope \cite{pp}, but our viewpoint and
interests are rather different from theirs. This family exhibits
several interesting features. It provides the first explicit examples
of a one-parameter family of PE metrics undergoing conic degeneration;
similarly, it also provides a simple nonproduct and nonhyperbolic
example (in fact, the first one known to us) of a family of PE metrics
`bending' along a codimension two submanifold. Both of these phenomena
occur for Einstein metrics with special holonomy, but these metrics
appear not to have any extra structure of this sort (even though the
building blocks do). Finally, we also obtain a one-parameter family of
PE metrics (albeit with interior edge singularities) for which the
conformal infinities are all equal to the standard conformal class on
the sphere; this stands in contrast to the fact that hyperbolic space
is known to be the unique smooth PE filling of this standard spherical
structure \cite{Q}.

Let us set this in a broader context. 
There is a natural asymptotic boundary problem for PE metrics: given the conformal infinity data $(Y,[\gamma_0])$, 
find a manifold $X$ with $\partial X = Y$ and a PE metric $g$ on $X$ with ${\mathfrak c}(g) = [\gamma_0]$. 
This has turned out to be quite difficult, and the current existence theory is mostly based on perturbation 
arguments, cf.\ \cite{An}, \cite{GL}, \cite{Bi}, \cite{MP}. The problem was originally posed by Fefferman 
and Graham \cite{FG} as a tool in conformal invariant theory; they showed the existence of formal (power 
series) solutions, which sufficed for their immediate purposes, but the global problem remains mostly
open. It is of interest not only in geometric analysis but also in string theory \cite{Bi2}. 

A closely related problem is to describe the moduli space of PE
metrics on a given manifold $X$; this is less constrained than the
analogous question in the compact setting, but just as in that case,
it leads naturally to the issue of compactness of families of PE
metrics. As usual, an indication for the failure of compactness is if
there exist families of PE metrics which develop singularities in the
limit, or more or less equivalently, if some family of PE metrics
$g_t$ on $X$ converges to an Einstein metric on a different space
$X_0$. A typical example is when $X$ is a complex manifold and the
limiting space is obtained by blowing down a divisor; the limiting
metric has a conic singularity at this singular point. Conic collapse
is known to be the only mode of degeneration for sequences of compact
four-dimensional Einstein manifolds under certain hypotheses
\cite{CT}, but in the PE setting there seems to be also the
possibility of cusp formation \cite{An}, although no examples of the
latter (for sequences of metrics on a fixed manifold) are
known. Sequences of Einstein metrics in higher dimensions may experience far worse types of 
degeneration, though it is reasonable to speculate that singularities of the limit 
should lie only in subsets of codimension four.

A key motivation for this project was provided by our efforts to
produce many examples of conically degenerating PE spaces using
analytic gluing techniques. Following many other successful
constructions of this type, such a construction should proceed by
starting with a PE manifold with isolated conic singularities,
replacing a neighbourhood of these singular points by scaled truncated
Ricci-flat ALE spaces, and then perturbing these ensembles to be
Einstein again. Although this scheme works easily if one assumes that
the component spaces -- in particular, the ALE one -- are
nondegenerate in the sense that there are no decaying Jacobi fields
for the linearized gauged Einstein operator, unfortunately there are
no known examples of nondegenerate Ricci-flat ALE spaces. A closer
examination indicates that this may be no accident, and that
degeneracy may indeed be a key feature which allows this conic
degeneration to occur. Our results in this direction are still
inconclusive, but this has made the examination of specific examples
all the more important. We hope that better understanding of the
spaces here will shed light on what to expect in general.

Now let us turn to the other singular behaviour in our special family of metrics.
Analytic constructions of families of Einstein metrics bending along a codimension
two submanifold seem to be even more obstructed, hence even harder; there is an extensive
literature on this in the three-dimensional hyperbolic setting, cf.\ \cite{HK}, \cite{W}
and \cite{MM}, but only very limited results have been obtained in higher dimensions \cite{Mo}. 
Hence once again the examples here may be of some interest in clearing the way toward a more
general construction. 

An interesting feature of this Page-Pope family of singular PE metrics is that there is
a nontrivial subfamily for which the conformal infinity (defined below) is the standard conformal 
structure on the sphere. This is slightly unexpected in light of the known uniqueness theorem when 
the interior of the PE space is smooth \cite{Q}, though presumably this happens here because of 
the negativity of the degree of the line bundle $L$ and the fact that the cone angle along
$Z$ is greater than $2\pi$ in these examples. The reason is that the proof in \cite{Q}
reduces the situation to an asymptotically flat one where the positive mass theorem can be
applied. This line of reasoning can be extended even when certain fairly mild interior singularities
are allowed, using \cite{Mi}, so the examples here show that there are genuine limits to
the possibility of such extensions. One is left with the problem of formulating the correct 
rigidity statement for PE metrics which are spherical at infinity but which have interior edge 
singularities;  this could have interesting consequences for extensions of the positive mass 
theorem too. 

It is necessary to introduce several definitions before stating our results precisely, and
this will be done in the next section. After that, we describe the ansatz for this family of
metrics and examine the various phenomena in different regions of the parameter space.

As noted already, this note is part of a larger effort of ours to
understand some of the analytic and geometric problems discussed
above; this project has been funded by the EPSRC grant
GR/S61522/01; R.M.\ was also supported by NSF grant DMS-0505709.

\section{Statements of results}

We begin by recalling the various classes of metrics with which we shall work, and then state
in detail the most interesting features of our family of metrics.

\subsection*{PE and edge metrics}
First, to expand on the discussion above, a {\em conformally compact}
space is the interior of a manifold with boundary $X$, endowed with a
metric of the form $g = \rho^{-2}\overline{g}$; here $\rho$ is a
defining function for $\del X$, so $\{\rho = 0\} = \del X$ and $d\rho
\neq 0$ there, and $\overline{g}$ is a Riemannian metric which is
smooth (or has some specified finite regularity) up to the
boundary. Any such metric is complete, and if $|d\rho|_{\overline{g}}
= 1$ at the boundary, then all sectional curvatures tend to $-1$ at
infinity.

A space $(X,g)$ is called Poincar\'e--Einstein (PE) if it is conformally compact and Einstein. It is 
usually convenient to normalize so that $\mbox{Ric}\,(g) = -(N-1)g$, $N = \dim X$. PE spaces are the 
asymptotically hyperbolic analogues of gravitational instantons, which include the Ricci flat ALE spaces 
we discuss briefly below. 

The conformal infinity of a PE metric $g$ (or indeed of any conformally compact metric) is the conformal class 
\[
{\mathfrak c}(g) = \left[\left.x^2 g\right|_{TY}\right]
\]
on $Y$; only the conformal class is well-defined since $x$ can be modified by any smooth nonvanishing 
factor. We say that $(X,g)$ bounds or fills $(Y,[h])$. 

Given any representative $\gamma_0 \in {\mathfrak c}(g)$, there is a unique associated special boundary 
defining function $x$ in some neighbourhood ${\mathcal U}$ of the boundary such that $|dx/x|^2_{g} \equiv 1$ 
there. Using the gradient flow of $x$, we identify ${\mathcal U}$ with $[0,\epsilon)_x \times Y$, and
can then write
\[
g = \frac{\rd x^2 + \gamma(x)}{x^2}
\]
where $\gamma(x)$ is a family of smooth metrics on $Y$ depending smoothly on $x$ with $\gamma(0) = \gamma_0$.
This is called the Graham-Lee (or standard) normal form for $g$. Note that this assumes that the Einstein
constant is normalized as above. 

Now we recall the definition of an edge metric. Let $X$ be a compact
smoothly stratified space with a dense top-dimensional smooth stratum,
and only one other lower dimensional stratum, $Z$, which is a compact
manifold disjoint from $Y$ and along which $X$ has an edge
singularity. Thus, by definition, some neighbourhood ${\mathcal V}$ of
$Z$ in $X$ is diffeomorphic to a bundle over $\pi: {\mathcal V} \to Z$
where each fibre is a truncated cone $C(F)$ over some compact smooth
manifold $F$. (Recall that $C(F)$ is the space obtained from
$[0,1)\times F$ by collapsing $\{0\}\times F$ to a single point $p$.)
Set $\dim X = N$ and $\dim Z = k$. (In practice below, $X$ will also
have a codimension one boundary disjoint from $Z$, but since this
definition is local near $Z$, that is of no concern.)

A Riemannian metric $g$ on $X$ is called an {\em edge metric} if in
the neighbourhood ${\mathcal V}$ it has the form
\[
g = g_0 + g', \qquad \mbox{where} \quad g_0 = \rd s^2 + s^2 h + \pi^* k.
\]
Here $s$ is the radial variable in the fibres of $L$ (with respect to
the given Hermitian metric), $h(s) \in {\mathcal C}^\infty([0,1]\times
F; S^2T^*F)$ is a smooth family of metrics on $F$, $k$ is a smooth
metric on $Z$, and $g'$ is a lower order term in the sense that
$|g'|_{g_0} \to 0$ (usually like some power of $s$) as $s \to
0$. Typically, various regularity assumptions are placed on $g'$.

The two cases which appear below are when $Z$ is a point, so that $X$
has isolated conic singularities (and a boundary), or when $Z$ has
codimension two, so that $h$ has the form $\alpha^2\, d\psi^2$, $\psi
\in S^1$ for some function $\alpha > 0$. For the metrics considered in
this paper, $\alpha$ is 
always constant along $Z$, and we refer to
the number $2\pi \alpha$ as the cone angle along this edge. In the
former case, in our examples, 
$X$ itself is diffeomorphic to a truncated cone $C(Y)$ and $F$ is
identified with the boundary $Y$.

The space $X$ can be identified topologically with a smooth manifold
with boundary where $Z$ is a smooth embedded submanifold if and only
if $F = S^{N-k-1}$; in this case, if an edge metric $g$ extends
smoothly over $Z$, then $h(0)$ is the standard metric on the sphere
(so $\alpha = 1$ when $k = N-2$). Of course, there are infinitely many
other compatibility conditions to ensure that the metric $g$,
expressed in Fermi coordinates around $Z$, is smooth on $X$. However,
when $g$ is Einstein, standard elliptic regularity arguments can be
used to prove that if $g$ is ${\mathcal C}^1$ near the submanifold
$Z$, then it is actually ${\mathcal C}^\infty$ across $Z$. We shall
use this principle below in several places without further comment.

Note that when $Y = S^{2n+1}$ with its standard metric, then $C(Y)$ is diffeomorphic to $\C^{n+1}$; more 
generally, if $Y = S^{2n+1}/\Z_k$, for some free linear action of $\Z_k$ on $S^{2n+1}$, then $C(Y)$ is the 
orbifold $\C^{n+1}/\Z_k$. 

\subsection*{Features of the Page-Pope metrics} 
The Page and Pope ansatz \cite{pp} requires the following data:
\begin{enumerate}
\item[(i)] A compact K\"ahler-Einstein manifold $(M,\hg)$ with $\dim_{\C}M = n$ and 
$c_1 > 0$; we set $\mathrm{Ric}(\hg) = \lambda \hg$.
\item[(ii)] A holomorphic line-bundle $L\to M$ equipped with Hermitian metric, with curvature 
form $-2\omega$, where $\omega$ is the K\"ahler form for $\hg$. Thus, if $Y$ denotes the associated 
$S^1$ bundle and $\theta$ is the connection $1$-form on $Y$, then $\rd\theta = -2\omega$. 
\end{enumerate}
We are following the conventions in \cite{pp} whereby $2[\omega] \in
H^2(M,2\pi i \Z)$, which then forces $(-2/\lambda) c_1(M) \in H^2(M,2
\pi i \Z)$ as well; this class corresponds to $c_1(L)$, of course, and gives
restrictions on the possible values of $\lambda$. In the key examples
below, $M = \C P^n$ and $\hg$ is some rescaling of the Fubini-Study
metric. Since $c_1(\C P^n) = (n+1)\omega$, we find that $L$ has degree
$k = -2(n+1)/\lambda$, so that $\lambda$ can only take on the values
$-(2n+2)/k$, where $k$ is an integer.  In any case, the value of
$\lambda$ fixes the topology of the total space of $L$, which we
denote $X$, and the boundary of its fibrewise radial compactification,
which is the $S^1$ bundle $Y$. We also denote by $Z\subset X$ the
zero-section and $X_0$ the blow-down of $X$ along $Z$. It is not hard
to verify that $X_0$ is naturally identified with the cone $C(Y)$.

We can now state the key features and properties of the families of
Einstein metrics we construct:

\begin{enumerate}
\item[(i)] Choose $L$ to be the canonical bundle $K\to M$. This corresponds to taking $\lambda=2$; 
  on the total space $X$ of $K$ there exists a family of smooth PE metrics $\{g_t\}_{t>0}$ such 
  that as $t \searrow 0$, $(X, g_t)$ converges to a PE metric on $C(Y)$ with an isolated
  conic singularity at the blow-down of $Z$.  The convergence is ${\mathcal C}^\infty$ away from 
  $Z$ and $Y$, but is smooth in a suitably modified sense even near these singular loci. There is 
  a limit as $t \searrow 0$ of the rescaled family $(X,t^{-1}g_t)$ which is a complete Ricci-flat 
  K\"ahler metric on $X$; the convergence is again smooth in an appropriately rescaled sense.
\item[(ii)] Suppose that $M = \C P^n$ and $\hg$ is the multiple of the Fubini-Study metric for which 
  $\lambda = 2n+2$. (This is the standard scaling in complex geometry.) Then $L$ is the degree $-1$ 
  bundle over $\C P^n$ and $X$ is biholomorphically equivalent to the complex blow-up of the origin 
  in $\C^{n+1}$; we also identify it with the interior of the compact manifold with boundary obtained by 
  taking the radial compactification of each fibre. Thus $Y = S^{2n+1}$ and $C(Y) = B$. There is a 
  family of PE metrics $\{g_t\}_{t>0}$ on this space $X$, each of which has an edge singularity along $Z$.
  As $t$ increases from $0$ to $\infty$, the cone angle along $Z$ increases monotonically from $2\pi(n+1)$ 
  to $\infty$. The diameter of $Z$ with respect to $g_t$ is of order $t$, so that as $t \searrow 0$, 
  $Z$ collapses to a point; in this limit, $g_t$ converges to the standard hyperbolic metric on $B$. 
  Notably, we may choose this family so that its conformal infinity $[\gamma_t]$ is the standard round conformal 
  structure on $S^{2n+1}$ for all $t$.
\item[(iii)] Now let $L$ be the line bundle of degree $-n-1$ over $\C P^n$; the Fubini--Study metric $\hg$ 
  is now normalized so that $\lambda =2$. The total space $X$ admits two distinct families of PE metrics, 
  $\{g_t\}$ and $\{g_t{\,'}\}$. The metrics $g_t$ constitute the family described in (i) above, so each is 
  smooth across $Z$ when $t>0$ and they converge as $t \searrow 0$ to a metric on $X_0$ with isolated conic
  singularity. The conformal infinity of $g_t$ varies with $t$, and in particular never equals the quotient 
  of the round conformal structure on $S^{2n+1}$.  On the other hand, the family $g_t{\,'}$ is the analogue of 
  the one in (ii): the conformal infinity of each $g_t{\,'}$ is the standard one on $S^{2n+1}/\Z_{n+1}$, and each
  metric has an edge singularity along $Z$. As $t$ increases, the cone angle along $Z$ increases monotonically
  from $2\pi (2n+1)/(2n+2)$ to infinity, and the diameter of $Z$ increases from $0$ to $\infty$ as well. In 
  particular, there is precisely one value of $t$ for which the cone angle is $2\pi$, so that the corresponding metric
  is smooth. As $t \searrow 0$, $(X,g_t{\,'})$ degenerates to the standard orbifold PE metric on $X_0 = B/\Z_{n+1}$, 
  i.e.\ the quotient of hyperbolic space by $\Z_{n+1}$. 
\end{enumerate} 
We also describe the Einstein metrics when the bundle $L \to \C P^n$ has any other negative degree too. \label{s2}

\section{The Page--Pope construction}

The following result is proved in \cite{pp}: 

\begin{prop}\label{ppp}
  Fix constants $c>0$, $\Lambda < 0$ and $r_1\geq 1$. Let $P(r)$ be
  determined by the equation
\begin{equation}\label{e2.18.6.5}
  \frac{\rd}{\rd r}(r^{-1}P(r)) = r^{-2}[|\Lambda|(r^2-1)^{n+1} +
  c^{-1}\lambda (r^2-1)^n], 
\end{equation}
and the condition $P(r_1) = 0$.  Let $(M,\hg)$, $\lambda$, $L$ and $X$
be as above, and suppose 
that $r$ is the radial variable in the fibres of $L$. Then the metric 
\begin{equation}\label{e1.17.6.5}
g = (r^2-1)^n P(r)^{-1} \rd r^2 + c^2 P(r)(r^2-1)^{-n}\theta^2+ c(r^2-1)\hg
\end{equation}
on $\{r> r_1\}\subset X$ is Einstein, with $\mathrm{Ric}\,(g) = \Lambda \, g$. Its rescaling
$\frac{|\Lambda|}{2n+1} g$ is Poincar\'e--Einstein as $r\to\infty$, with conformal infinity
\begin{equation}\label{e1.31.8.7}
\left[\frac{c|\Lambda|}{2n+1}\theta^2 + \hg\right].
\end{equation}
If $r_1>1$, then $g$ has an edge singularity along the zero section $Z$ in $X$, and in fact 
is asymptotic near this singular locus to a constant multiple of
\begin{equation}\label{e2.31.8.7}
\rd s^2  + \alpha^2s^2\theta^2 +  \beta^2\,  \hg,
\end{equation}
where $s$ is the radial distance function. Here
\begin{equation}\label{e3.31.8.7}
\alpha = \frac{c|\Lambda|}{2r_1}(r_1^2-1) + \frac{\lambda}{2r_1},\;\;
\beta^2 = \frac{\alpha^2 (r_1^2-1)}{2}.
\end{equation}
On the other hand, if $r_1=1$, then $\hg$ defines a PE metric on $X_0$ with a conic singularity 
at the blowdown of $Z$, near which it is asymptotic to the model
\begin{equation}\label{e4.31.8.7}
\rd s^2 + s^2\left(\frac{\lambda^2}{(2n+2)^2}\theta^2 + \frac{c\lambda}{2n+2} \hg \right).
\end{equation}
\end{prop}

We make some comments. In order that \eqref{e1.17.6.5} define a metric when $r > r_1$, it is obviously
necessary that $r>1$ and $P(r)>0$ in this range. This explains the restriction $r_1 \geq 1$; furthermore,
the RHS of \eqref{e1.17.6.5} is positive then, so the positivity of $P(r)$ when $r > r_1$ is immediate
from the ODE and boundary condition. 

The proof that $g$ is Einstein is a straightforward computation, cf.\ \cite{pp}. Although their notation has 
been followed closely, we replaced $P(r)$ by $(-1)^nP(r)$ and $c$ by $-c$ in order to make
certain quantities positive in what follows.

The limit $r\to\infty$ corresponds to the conformal infinity of $g$. To see this, note that for large $r$, the RHS of
\eqref{e1.17.6.5} equals $|\Lambda|r^{2n} + O(r^{2n-1})$, hence 
\begin{equation}\label{e2.28.8.7}
g \sim \frac{2n+1}{|\Lambda|}\frac{\rd r^2}{r^2} + c^2\frac{|\Lambda|}{2n+1}r^2\theta^2 + cr^2\hg, \qquad r \gg 0.
\end{equation}
Setting $x=r^{-1}$ and rescaling $g$ by the appropriate constant gives a PE metric with conformal infinity 
$[c|\Lambda|/(2n+1)\theta^2 + \hg]$.

To verify the other statements, we must analyze $P(r)$ near $r=r_1$. If $r_1>1$ and we write 
$r = r_1 +s^2$, then for $s \ll 1$, 
$$
g \sim 4\frac{(r_1^2-1)^n}{P'(r_1)}\left(\rd s^2 + \left(\frac{cP'(r_1)}{2(r_1^2-1)^n}\right)^2s^2\theta^2
+ \frac{cP'(r_1)}{4(r_1^2-1)^{n-1}} \hg\right).
$$
The ODE gives
$$
P'(r_1) = \frac{1}{r_1}(|\Lambda|(r_1^2-1)^{n+1}  + c^{-1}\lambda(r_1^2-1)^{n})>0,
$$
so $g$ has an edge singularity along $Z$. Its cone angle is $2\pi$ times
\[
\alpha = \frac{c P'(r_1)}{2(r_1^2 - 1)^n} = \frac{c|\Lambda|}{2} r_1 + \frac{\lambda - c |\Lambda|}{2} \frac{1}{r_1},
\]
as claimed. For certain values of $c$, $\Lambda$ and $\lambda$, this
need not be monotone in $r_1$, but in any 
case it certainly satisfies
\[
\lim_{r_1 \searrow 1} \alpha = \frac{\lambda}{2}, \qquad \lim_{r_1
  \nearrow \infty} \alpha = \infty
\]
for fixed $(c,\Lambda, \lambda)$.

When $r_1=1$, then for small $s$,
$$
P(1 + s^2) \sim \frac{2^{n+1}s^{2n+2}}{n+1},
$$
and the claim about the regularity of the metric follows by substituting this into \eqref{e1.17.6.5}.

This construction involves four parameters, $\lambda$, $c$, $\Lambda$ and $r_1$. Amongst these, $\lambda$ 
should be regarded as fixed in advance, since it determines the topology, but it is convenient to maintain 
flexibility with the others. Fixing $\lambda$ and $r_1$, the scaling
$$
(c,\Lambda) \mapsto (ac, a^{-1}\Lambda)
$$
has the effect of replacing $P$ by $aP$, which in turn changes $g$ to
$a^{-1}g$. Using this freedom, we can always assume that $\Lambda = -2n-1$, which is the standard Riemannian
normalization for Einstein metrics, and puts $g$ into standard form as $r\to\infty$.  However,
this still leaves two effective free parameters.

It is sometimes more convenient to fix the conformal infinity of $g$,
which corresponds to fixing $c$.  In this case, the only remaining
free parameter is $r_1$. On the other hand, demanding that $g$ is
smooth across $Z$ for any given $r_1>1$ fixes $c$ in a possibly
different way.  Hence we cannot expect to prescribe both the conformal
infinity {\em and} also require that the metric be smooth across
$Z$. We describe this in more detail in the next section.

\section{A degenerating family of smooth PE metrics}

According to the final paragraph of the last section, after fixing
$\lambda$ and setting $\Lambda=-2n-1$,  
there is a two-parameter family of PE metrics parametrized by $r_1$ and $c$. From now on, let
$r_1 = 1 + t$ and write $P_t(r)$ and $g_t$ for the corresponding solution and metric.

As explained earlier, in order to understand when $g_t$ is smooth across $Z$, it suffices to 
consider the limiting metric on the normal circles, and in particular to check that the coefficient 
$\alpha$ in \eqref{e3.31.8.7} is equal to $1$. This is the condition
\begin{equation}\label{e1.1.9.7}
c|\Lambda| = \frac{2r_1-\lambda}{r_1^2-1} \Rightarrow c = c_t = \frac{1 + t - \lambda/2}{(2+t)(2n+1)},
\end{equation}
which implies in particular that $r_1 \geq \max(1,\lambda/2)$. Hence, for any fixed $\lambda$, there is 
exactly one smooth metric in this family for each admissible value of $t$. In order to continue this
family of smooth metrics to the conical limit at $t=0$, i.e.\ $r_1 = 1$, we must have $\lambda \leq 2$. 
In other words, when $\lambda \leq 2$, there is a family of smooth PE metrics $g_t$, $t > 0$, such that 
as $t \searrow 0$, $g_t$ converges to a PE metric on $X_0 = C(Y)$ with isolated conic singularity. 

Let us consider the case $\lambda = 2$ in more detail. We shall analyze the behaviour of this degenerating 
family more closely by a rescaling of the space, which we understand by rescaling the radial coordinate. 
Since the diameter of $Z$ with respect to $g_t$ is of order $t$, it is natural to consider the limit of 
the rescaled family $t^{-1}g_t$. (This is the same as the rescaling considered in \S3 of \cite{pp}, where the 
authors allow $c \to \infty$.) Set $\rho^2 = c(r^2-1)$ and $U(\rho) = cP_t(r)/(r^2-1)^{n+1}$. Now apply the scaling 
\[
(c_t,\Lambda) = (c_t,-(2n+1)) \mapsto (c_t/t, -(2n+1)t)
\]
and let $t\to 0$.  In terms of the variable $\rho$ and the function $U$, the limiting metric
is equal to 
\begin{equation}\label{e1.29.8.7}
g_\infty = \lim_{t\to 0} t^{-1}g_t = U(\rho)^{-1}\rd \rho^2 + U(\rho)\rho^2\theta^2 + \rho^2 \hg,
\end{equation}
where 
\begin{equation}\label{e5.18.6.5}
\frac{\rd}{\rd \rho}(\rho^{2n+2}U)  = 2 \rho^{2n+1}, \qquad U(2/\sqrt{2n+1}) = 0.
\end{equation}
Note that the lower limit for $\rho$ is $\rho_1 = 2/\sqrt{2n+1}$ since this is the limiting solution
to $\rho_1^2 = (c_t/t)((1+t)^2 - 1)$ as $t \to 0$. As at the end of \S3 of \cite{pp}, $g_\infty$ is a 
Ricci-flat K\"ahler metric on $X$. 

We must verify that $g_\infty$ is smooth at $\rho = \rho_1=2/\sqrt{2n+1}$. The solution to (\ref{e5.18.6.5}) is
\begin{equation}\label{e11.29.8.7}
U(\rho) = \frac{1}{n+1}(1 - (\rho_1/\rho)^{2n+2}),
\end{equation}
so with $\rho = \rho_1 + s^2$, the metric becomes 
$$
\rho_1\left(2\rd s^2 + 2 s^2 \theta^2 + \rho_1 \hg\right)
$$
for small $s$, which is regular at $s=0$.

This fills in the details of the construction of the first family of
metrics described at the end of \S\ref{s2}.

\section{PE metrics on line-bundles over $\C P^n$}

Now we turn to the second family of metrics from the end of
\S\ref{s2}.

Let
$(M,\hg)$ be $\C P^n$ with Fubini--Study metric normalized by
$\mbox{Ric}\,(\hg) = \lambda \hg$; as remarked earlier, the degree of $L$
equals $-(2n+2)/\lambda$, so $\lambda = (2n+2)/k$ for some $k \in
{\mathbb N}$.

Both the flat metric on $\R^{2n+2}$ and the hyperbolic metric on the
ball $B^{2n+2}$ can be recovered from this construction. Indeed, to
get the topology right, we must choose $k = 1$, so $\lambda = 2n+2$.
The flat metric is Ricci-flat and K\"ahler, so it arises from the same
sort of limit of $g_t/t$ as above, and corresponds to the constant
solution $U(\rho) = 1$. The metric is
\begin{equation}\label{e3.2.9.7}
g_{{\mathrm{flat}}} = \rd \rho^2 + \rho^2\theta^2 + \rho^2 \hg,\; \qquad \mbox{where}
\qquad \mbox{Ric}\,(\hg) = (2n+2)\hg.
\end{equation}
Similarly the hyperbolic metric arises as the conical limit with $\Lambda = -2n-1$, $c=1$
and $r_1 \searrow 1$. Notice that by varying $c$, we get a family of conical PE metrics
filling the conformal classes $[c\theta^2 + \hg]$ on $S^{2n+1}$; the standard round metric 
appears in this family if and only if $c=1$.

On the other hand, we can make a deformation starting from the
standard hyperbolic metric holding the boundary conformal structure
fixed by varying $r_1$. By Proposition~\ref{ppp}, this gives a metric
on $X$ with an edge singularity along the zero set $Z$. It is not hard
to check that if $\Lambda=-2n-1$ and $c=1$, then the quantity $\alpha$
of \eqref{e3.31.8.7} never equals $1$, and hence $g$ is not smooth
across $Z$.

For the bundle of degree $-k$ over $\C P^n$, the discussion above goes through with only minor modifications.  We fix 
$\Lambda =-2n-1$ and $\lambda = (2n+2)/k$, and take $c= 1/k$ so that the conformal infinity (\ref{e1.31.8.7}) is 
always equal to the round conformal structure on $S^{2n+1}/\Z_k$.  As $r_1$ varies between $1$ and $\infty$,
the cone angle increases monotonically from $2\pi (n+1)/k$ to $\infty$, and in the limit as $r_1 \searrow 1$, 
we get the standard orbifold hyperbolic metric on $B/\Z_k$. The ALE rescaling (or `near horizon' limit)
is a Ricci-flat K\"ahler space with edge metric along $Z$. 
As before, we may deform this by increasing $r_1$, holding the other quantities fixed, and again the 
result is a family of PE metrics with fixed conformal infinity and with an edge singularity along the zero-section 
of the line bundle.

Perhaps the most interesting case is $k=n+1$ where this family of metrics coexists with the degenerating family of 
smooth PE metrics considered in the previous section.  These two families have conformal infinities which are quite 
different from each other.  This provides some circumstantial evidence that a small perturbation of the round
conformal structure on $S^{2n+1}/\Z_{n+1}$ does not bound a smooth PE metric on the total space of the line bundle 
$L$ of degree $-n-1$ over $\C P^n$.  This completes our discussion of the examples described at the end of \S\ref{s2}.

\section{Further directions}
The explicit Poincar\'e-Einstein spaces here, both singular and nonsingular, are interesting 
in their own right, of course, but may be particularly useful in providing 
intuition for some analytic constructions which would show that conic
degeneration and bending are robust phenomena. As noted in the
introduction, based on heuristic arguments, it appears that a direct
gluing construction to construct conically degenerating families of PE
metrics is always obstructed by the existence of a cokernel for the
relevant elliptic operator. If this is true, it would indicate that
the limiting conic space would have to be degenerate, i.e.\ admit $L^2$
Jacobi fields for the linearized gauged Einstein operator, as described in \cite{MP}.  
Thus it would be very interesting to see whether the limiting conic PE spaces obtained 
here are degenerate in this sense.

Another question is whether it is possible to take a smooth Einstein metric and start bending it along 
a codimension two submanifold in a family of Einstein metrics with edge singularities. In some of the examples 
here, this actually occurs. It would be very interesting to understand the correct obstructions for
this problem.


\begin{thebibliography}{99}
\bibitem{An} M.T. Anderson {\em Geometric aspects of the AdS/CFT correspondence} in {\sl AdS/CFT correspondence: 
Einstein metrics and their conformal boundaries} IRMA Lect. Math. Theor. Phys., 8, Eur. Math. Soc., 
Zürich (2005), 1--31. 
\bibitem{Bi} O. Biquard {\em Asymptotically symmetric Einstein metrics} Translated from the 2000 French original 
by Stephen S. Wilson. SMF/AMS Texts and Monographs, 13. AMS, Providence, RI; 
Société Mathématique de France, Paris, 2006.
\bibitem{Bi2} O. Biquard, ed. {\em  AdS/CFT correspondence: Einstein metrics and their conformal boundaries} 
IRMA Lect. Math. Theor. Phys., 8, Eur. Math. Soc., Zürich, 2005.
\bibitem{CT} J. Cheeger and G. Tian, {\em Curvature and injectivity radius estimates for Einstein 4-manifolds}  
J. Amer. Math. Soc.  19 no. 2 (2006) 487--525.
\bibitem{CY} S.Y. A. Chang and P. Yang {\em Boundary regularity of Bach-flat metrics}, in preparation.
\bibitem{FG} C. Fefferman and C.R. Graham {\em Conformal invariants} in {\sl The mathematical heritage of 
Élie Cartan} Astérisque (1985),  Numero Hors Serie, 95--116. 
\bibitem{GL} C.R. Graham and J. Lee {\em Einstein metrics with prescribed conformal infinity on the ball} 
Adv. Math. {\bf 87} no. 2 (1991) 186--225. 
\bibitem{HK} C. Hodgson and S. Kerckhoff {\em Rigidity of hyperbolic cone-manifolds and hyperbolic Dehn surgery}
 J. Diff. Geom. {\bf 48} no. 1 (1998), 1--59.
\bibitem{MP} R. Mazzeo and F. Pacard {\em Maskit combinations of Poincaré-Einstein metrics}  
Adv. Math.  {\bf 204} no. 2 (2006), 379--412.
\bibitem{MM} R. Mazzeo and G. Montcouquiol {\em The infinitesimal Stoker conjecture and deformations of
cone manifolds}, In preparation.
\bibitem{Mi} P. Miao {\em Positive mass theorem on manifolds admitting corners along a hypersurface} 
Adv. Theor. Math. Phys. {\bf 6} no. 6 (2002), 1163--1182.
\bibitem{Mo} G. Montcouquiol {\em D\'eformations de m\'etriques Einstein sure les vari\'et\'es \`a 
singularit\'e conique}. Thesis, Universit\'e Paul Sabatier -- Toulouse III, 2006.
\bibitem{pp} D.\ Page and C.\ Pope {\em Inhmogeneous Einstein metrics on complex line  bundles},
Class. Quantum Grav., {\bf 4} (1987), 213--225. 
\bibitem{Q} J. Qing {\em On the rigidity for conformally compact Einstein manifolds} Int. Math. Res. Not. 
no. 21 (2003) 1141--1153.
\bibitem{W} H. Weiss {\em Local rigidity of 3-dimensional cone-manifolds} J. Diff. Geom. {\bf 71}
no. 3 (2005) 437--506.
\end{thebibliography}
\end{document}